 \documentclass[pmlr,twocolumn,10pt]{jmlr} 

 \usepackage{mathtools} 

\usepackage{booktabs}
\usepackage{siunitx}

\usepackage{pgfplots}
\pgfplotsset{compat=1.18} 
\usepackage{tikz}

\DeclareMathOperator*{\argmax}{arg\,max}

\theorembodyfont{\upshape}
\theoremheaderfont{\scshape}
\theorempostheader{:}
\theoremsep{\newline}

\jmlrvolume{LEAVE UNSET}
\jmlryear{2023}
\jmlrsubmitted{LEAVE UNSET}
\jmlrpublished{LEAVE UNSET}
\jmlrworkshop{Machine Learning for Health (ML4H) 2023} 

 \title[Learning Control Policies of HH Neuronal Dynamics]{Learning Control Policies of Hodgkin-Huxley Neuronal Dynamics}

 \author{
  \Name{Malvern Madondo} \Email{malvern.madondo@emory.edu}\\
  \addr Department of Computer Science
  \AND
  \Name{Deepanshu Verma} \Email{deepanshu.verma@emory.edu}\\
  \Name{Lars Ruthotto} \Email{lruthotto@emory.edu}\\
  \addr Department of Mathematics
  \AND
  \Name{Nicholas {Au Yong}} \Email{nicholas.au.yong@emory.edu}\\
  \addr Department of Neurosurgery
  \AND
  \addr Emory University, Atlanta, GA, USA
 }
 
\begin{document}

\maketitle

\begin{abstract}
We present a neural network approach for closed-loop deep brain stimulation (DBS). We cast the problem of finding an optimal neurostimulation strategy as a control problem. In this setting, control policies aim to optimize therapeutic outcomes by tailoring the parameters of a DBS system, typically via electrical stimulation, in real time based on the patient's ongoing neuronal activity. We approximate the value function offline using a neural network to enable generating controls (stimuli) in real time via the feedback form. The neuronal activity is characterized by a nonlinear, stiff system of differential equations as dictated by the Hodgkin-Huxley model. Our training process leverages the relationship between Pontryagin's maximum principle and Hamilton-Jacobi-Bellman equations to update the value function estimates simultaneously. Our numerical experiments illustrate the accuracy of our approach for out-of-distribution samples and the robustness to moderate shocks and disturbances in the system.
\end{abstract}


\section{Introduction}
\label{sec:intro}

Disruptions in the normal functioning of neurons are an early symptom of many neurological conditions such as Parkinson's Disease \citep{brown2003oscillatory, wichmann2016deep}. Controlling neuronal dynamics has potentially high practical importance in developing neurotechnology and improving health-related quality of life. In neuromodulatory interventions such as deep brain stimulation, physicians aim to alter neuronal activity by delivering electrical pulses to target locations in the brain via a battery-powered pulse generator or neurostimulator \citep{arlotti2016adaptive, koeglsperger2019deep}.

A key clinical objective is to learn an effective stimulus input with optimal therapeutic benefits for the patient while reducing unwanted side effects. Traditionally, physicians have relied on ad hoc approaches to learn the optimal stimulation strategy over multiple clinical visits \citep{mohammed2018toward, parastarfeizabadi2017advances, yu2020review}. However, this approach has limitations as the neurostimulator continuously delivers a fixed pattern of stimulation, reducing the battery life in the long run and potentially being unnecessary during normal neuronal function \citep{shirvalkar2018closed, meidahl2017adaptive}. Despite recent advancements in neurotechnology, the development of adaptive stimulation strategies remains a challenging research area that requires real-time control algorithms based on the brain's current clinical state \citep{carron2013closing, parastarfeizabadi2017advances, wang2007bifurcation}.

In this work, we tackle the biomedical problem of learning a stimulation strategy by formulating it as an optimal control problem. This framework paves the way for developing adaptive closed-loop neuromodulation strategies and is promising for advancing neurostimulation therapies. We characterize neuronal activity with the Hodgkin-Huxley model \citep{hodgkin1952quantitative}, a biophysically accurate nonlinear system of differential equations that captures the essential dynamics of neuronal behavior. By combining Pontryagin's maximum principle (PMP) and Hamilton-Jacobi-Bellman (HJB) equation, we establish a solid mathematical foundation that allows us to obtain approximately optimal controls in real time. While the PMP is a local solution method and is sensitive to the initial states and perturbations in the optimal trajectories, it provides first-order necessary conditions for optimality. In contrast, the HJB equation yields controls in feedback form and provides sufficient conditions for optimality. However, solving the HJB equations in higher dimensions is notoriously challenging due to the curse of dimensionality (CoD). To overcome these limitations, we adopt a semi-global neural network (NN) approach that approximates the value function and learns a policy across a diverse set of initial conditions and the state space likely to be encountered when following optimal trajectories. This approach has been successfully demonstrated in various control problems in economics \citep{han2018solving}, trajectory planning \citep{onken2021multi}, stochastic optimal control \citep{pereira20a, li2022neural, Exarchos-2018}, and other domains.

\paragraph{Contributions} Our main contributions are as follows:
\begin{itemize}
    \item We provide a rigorous and comprehensive control formulation for designing an optimal neurostimulation strategy as a control problem. This approach enables the application of control theory and machine learning to address complex biological systems, leading to innovative approaches in both scientific research and medical practice.
    \item We establish a concrete link between the learning problem and optimal control theory, specified by the relation between PMP and HJB equations. We approximate the value function using a neural network satisfying the HJB, from which the optimal (stimulation) control can be recovered in real time via the feedback form. This approach is useful because it is transferable to closed-loop control problems in domains beyond healthcare. 
\end{itemize}

\section{Neuronal Dynamics}
\label{sec:background}

Developed in the 1950s, the Nobel-prize-winning Hodgkin-Huxley (HH) equations are widely considered the gold standard model of neuronal dynamics in computational neuroscience \citep{hodgkin1952quantitative}. The HH model ascribes the action potential (spike) generation in a neuron, i.e., the moment of maximum membrane potential/voltage, to the movement of fast depolarizing and slow hyperpolarizing ionic currents.

\begin{figure}[htbp]
\floatconts
  {fig:normal_HH_dynamics}
  {\caption{Evolution of the membrane potential and gating variables of \textbf{normal} HH equations, with no controls/stimuli, i.e., $u = 0, \, \forall \, t > 0$, and initial state $\vec{x} = [0, 0, 0, 0]^\top$.}}
  {%
   \subfigure[Membrane Potential][c]{\label{fig:normal_HH_action_potential}%
        \includegraphics[width=0.7\linewidth]{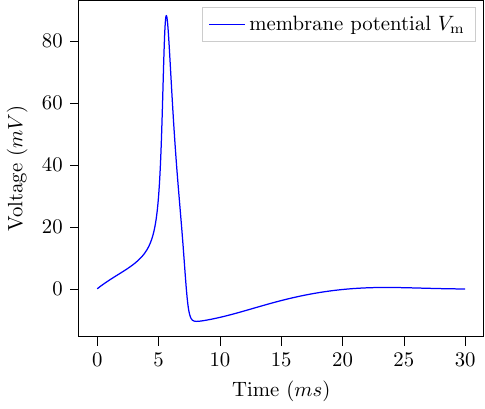}
    }
    \qquad
   \subfigure[Gating Variables][c]{\label{fig:normal_HH_gating_variables}%
        \includegraphics[width=0.7\linewidth]{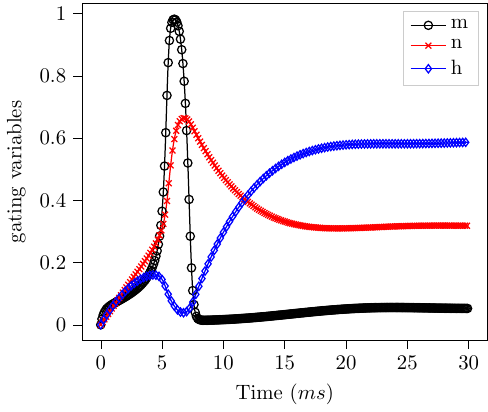}
    }
  }
\end{figure}

We define the state dynamics using the HH model, which describes the electrophysiological activity of neurons and accounts for currents flowing through the neuronal membrane across different ion channels, including sodium, potassium, and leakage channels. Leak channels comprise all other ions with slower dynamics, such as chloride ions. The HH model is characterized by the following first-order dynamical system, with state variables comprising of a neuron's  membrane potential (\figureref{fig:normal_HH_action_potential}), $V_{\rm m}$, and three gating variables (\figureref{fig:normal_HH_gating_variables}) representing probabilities of sodium activation, $m$, potassium activation $n$, and sodium inactivation, $h$

\begin{align}\label{eqn:dynamics}
\begin{split}
\frac{d \vec{z}}{dt}(t) &= f\big(t, \vec{z}(t)\big) + \vec{e}_1 u(t) \\
\vec{z}(0) &= \vec{x}, \quad \vec{e}_1 = [1, 0, 0, 0]^\top, \quad 0 \leq t \leq T
\end{split}
\end{align}
where $\vec{z}(t) = \left[V_{\rm m}, m, n, h\right]^\top \in \mathbb{R}^4$ denotes the state variable, with $\vec{x}$ being the initial state of the system. Here, $T$ denotes the fixed final time horizon, and the control variable $u(t): [0, T] \to \mathbb{R}$ represents the external current/stimulus provided as input by a controller at time $t$. The function $f: [0,T] \times \mathbb{R}^4 \to \mathbb{R}^4$ describes the evolution of state dynamics of the HH model and can be written as 
\begin{align}\label{eqn:f}
    \begin{multlined}
        f\left(t, \vec{z}(t)\right) =\\
        \left(\begin{array}{l}
            f_0(t,\vec{z})\\
            \alpha_m\big(z_0\big)(1 - z_1) - \beta_m\big(z_0\big)z_1\\
            \alpha_n\big(z_0\big)(1 - z_2) - \beta_n\big(z_0\big)z_2\\
            \alpha_h\big(z_0\big)(1 - z_3) - \beta_h\big(z_0\big)z_3
        \end{array}\right),
        \end{multlined}
\end{align}
where
\begin{align*}
    f_0\left(t, \vec{z}\right) = -\frac{1}{C_{\rm m}} \left( \begin{multlined} \Bar{g}_{\rm Na}z_1^3z_3(z_0 - E_{\rm Na}) \\ + \Bar{g}_{\rm K}z_2^4(z_0 - E_{\rm K}) \\ + \Bar{g}_{\rm l}(z_0 - E_{\rm l}) \end{multlined} \right).
\end{align*}
 Here, $C_{\rm m}$ is the membrane capacitance and $\Bar{g}_{\rm Na}, \Bar{g}_{\rm K}, \Bar{g}_{\rm l}$ represent the maximum conductance of the sodium and potassium and the leak currents, respectively. The variables $E_{\rm Na}, E_{\rm K}, \text{ and } E_{\rm l}$ represent the equilibrium potentials of the corresponding ions. The functions $\alpha_j$ and $\beta_j$, $j \in \{1, 2, 3\}$ defined in \tableref{table:hh_gating_variables}, are voltage-dependent rate functions for each of the gating variables. The default values and units of all the parameters are obtained from \citet{hodgkin1952quantitative} and listed under the ``normal" neuronal activity category in \tableref{table:HH_params}. \figureref{fig:normal_HH_dynamics} shows the state trajectory of typical HH neuronal activity.

The HH model (\equationref{eqn:dynamics}) is known to be stiff since $V_{\rm m}$ and $m$ change rapidly relative to $n$ and $h$. The rapid depolarization (membrane potential becomes more positive) and subsequent repolarization (membrane potential becomes more negative) of the neuron's membrane potential results in the spike/peak (formally known as the action potential) in \figureref{fig:normal_HH_dynamics}. On the one hand, the depolarization phase of the spike is largely driven by the influx of sodium ions due to the activation of the sodium channels. On the other hand, the repolarization phase is driven by the subsequent activation of the potassium channels and inactivation of the sodium channels. The general flow of ions, such as sodium and potassium, across the neuron’s membrane, plays a crucial role in transmitting electrical signals in neurons. We refer to \citet{gerstner2002spiking, koch2004biophysics, miller2018introductory} for detailed explanations of neuronal dynamics. 
\begin{table}[hbtp]
\floatconts
    {table:hh_gating_variables}
    {\caption{Parameters of the gating variables in the HH equations. The parameters $\alpha$ and $\beta$ depend only on the voltage, $z_0$.}}
    {\begin{tabular}{c c c}
    \toprule
    $j$ &  $\alpha_j$ &  $\beta_j$ \\ 
    \midrule
    $z_1$ & $\displaystyle \frac{2.5 - 0.1 z_0}{\exp(2.5 - 0.1 z_0) - 1}$ & $4 \exp(-\frac{z_0}{18})$\\
    \abovestrut{4.2ex}
    $z_2$ & $\displaystyle \frac{0.1 - 0.01 z_0}{\exp(1-0.1 z_0) - 1}$ & $0.125 \exp(-\frac{z_0}{80})$\\
    $z_3$ & $0.07 \exp(-\frac{ z_0}{20})$ & $\displaystyle \frac{1}{\exp(3 - 0.1 z_0) + 1}$ \\
    \bottomrule
    \end{tabular}}
\end{table}

\begin{figure*}[htbp]
\floatconts
{fig:pathological_HH_dynamics}
{\caption{Evolution of the membrane potential and gating variables of \textbf{pathological} HH equations, with no controls/stimuli, i.e., $u = 0, \, \forall \, t > 0$, and initial state $\vec{x} = [0, 0, 0, 0]^\top$.}} 
{%
    \subfigure[Membrane Potential][c]{\label{fig:pathological_HH_action_potential}%
        \includegraphics[width=0.35\linewidth]{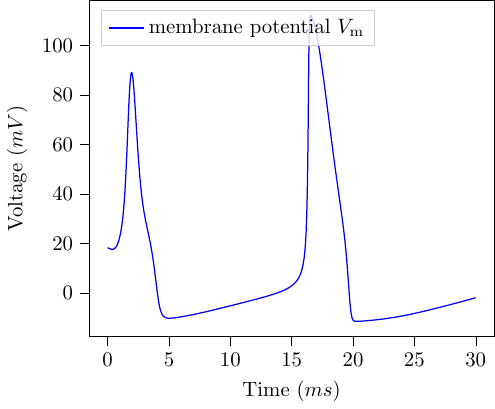}
    }
   \subfigure[Gating Variables][c]{\label{fig:pathological_HH_gating_variables}%
        \includegraphics[width=0.35\linewidth]{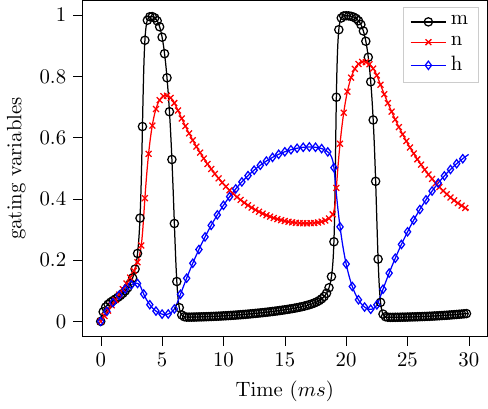}
    }
}
\end{figure*}

\begin{table}[hbtp]
\floatconts
{table:HH_params}
{\caption{Nominal parameter values of the HH model under normal and pathological conditions.}}
{
    \begin{sc}
    \begin{tabular}{lcccr}
    \toprule
    Parameter & Normal & Pathological \\
    \midrule
    $Cm$    & 1.0 & 1.0\\
    $\Bar{g}_{\rm Na}$ & 120.0 & 380.0\\
    $\Bar{g}_{\rm K}$    & 36 & 36.0\\
    $\Bar{g}_{\rm l}$     & 0.3 & 0.3\\
    $E_{\rm Na}$     & 115.0 & 115\\
    $E_{\rm K}$      & -12.0 & -12.0\\
    $E_{\rm l}$     & 10.613 & 10.613\\
    \bottomrule
    \end{tabular}
    \end{sc}
}
\end{table}

\subsection{Computational Approaches to Neurostimulation}
\label{subsec:comp_models}

Interest in applying machine learning and optimal control to neuroscience has soared recently. While the machine learning community generally favors a data-driven approach \citep{boutet2021predicting, Lu2020application, narayanan2023hh, peralta2021machine}, a model-based control approach offers opportunities to not only study the effect of external stimuli on neuronal dynamics but to also understand the computational properties of neural systems further \citep{breakspear2017dynamic, catterall2012hodgkin, salfenmoser2022nonlinear}.

Recent studies on developing computational models for neurostimulation generally focus on the dynamics of a single-compartment neuron model, or its extension, the multi-compartment model, which accounts for the spatial dimensionality of neurons and provides a more detailed representation of neuronal behavior \citep{deco2008dynamic, frohlich2005feedback, so2012relative, van2009mean}. The underlying dynamics in both compartment models are often described by HH-type equations (\equationref{eqn:f}). In this work, we focus on the single-compartment HH model describing the neuronal activity of a single neuron. We utilize the standard parameters of the HH model to represent normal neuronal dynamics. 

\paragraph{Pathological Activity} Additionally, we introduce a pathological state by randomly distorting one or more parameters, such as sodium conductance $\Bar{g}_{\rm Na}$, as outlined in \tableref{table:HH_params}. This has the effect of generating abnormal spikes in the membrane potential (see \figureref{fig:pathological_HH_action_potential}) and gating variables (\figureref{fig:pathological_HH_gating_variables}). This approach draws inspiration from a recent study by \citet{narayanan2023hh} that proposed a machine learning-based deep brain stimulator for controlling epileptic seizures. While our application differs from that study, the idea of simulating pathological conditions by manipulating the parameters of the HH model aligns with our objectives. \figureref{fig:normal_HH_dynamics} shows the evolution of the HH model under pathological conditions.

\section{Optimal Control of Neuronal Dynamics}
\label{sec:formulation}

In this section, we target the control of pathological neuronal activity (\figureref{fig:pathological_HH_dynamics}) by determining a control policy that can effectively restore normal functioning. We adopt an optimal control formulation of the neurostimulation problem that aims at restoring neuronal function (\figureref{fig:normal_HH_dynamics}) while expending less energy. Recent research enables this by combining function-approximating powers of NN and OC theory to solve high-dimensional control problems using neural networks \citep{onken2021neural, kunisch2021semiglobal}.

\subsection{Problem Setting}
For a fixed initial state $x$ and a finite time-horizon $T$, consider the optimal control problem constrained by the HH-based nonlinear system dynamics in \equationref{eqn:dynamics}. 
\begin{figure}[htbp]
    \floatconts
    {fig:normal_vs_pathological_action_potential}
    {\caption{Membrane potentials, $V_{\rm m}$, for normal condition ($\Bar{g}_{\rm Na} = 120$) vs a pathological condition ($\Bar{g}_{\rm Na} = 380$).}}
    {
    \includegraphics[width=0.7\linewidth]{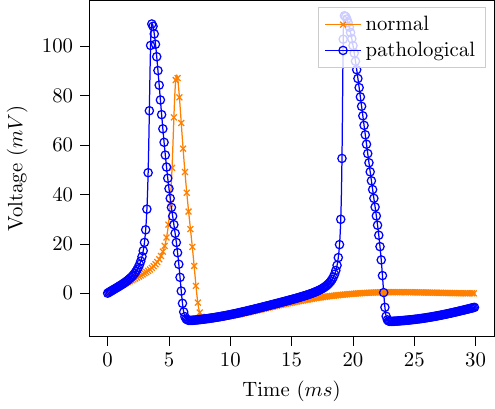}
    }
\end{figure}

We aim to find an optimal control (stimulus) $u^*: [0,T] \to \mathbb{R}$ that drives the system toward a reference/target state, $\vec{z}^*$, while incurring the minimum possible cost. To this end, we define the control objective functional (cost) $J$ as
\begin{align}\label{eqn:cost_function}
    J(t, \vec{z}, u) &= G\big(\vec{z}(T)\big) + \int_t^T L\big(s, \vec{z}(s), u(s)\big) ds,
\end{align}
with the terminal cost $G: \mathbb{R}^4 \to \mathbb{R}$ defined by
\begin{align}
    G(z) &= \frac{1}{2} \|\vec{z} - \vec{z}^*(T)\|^2,
\end{align}
and the Lagrangian, or running cost, $L: [0, T] \times \mathbb{R}^4 \times \mathbb{R} \to \mathbb{R}$ given by
\begin{align}
    L(t, \vec{z}, u) &= \lambda \|u\|^2 + Q \frac{1}{2} \|\vec{z} - \vec{z}^*(t)\|^2.
\end{align}

While the terminal cost, $G$, penalizes the distance between the final state $\vec{z}(T)$ and the given target terminal state of the system $\vec{z}^*(T)$, $L$ accumulates the cost of controlling the system and expending energy at each time step. Here, $Q$ penalizes the tracking term and is set to $200$ in our experiments. $\lambda$ is a problem-specific constant for electrode impedance used to limit the energy of the control signal. It is typically set to $0.5$ \citep{fleming2020simulation}.

Overall, we seek to minimize $J$ over all admissible controls $u \in \mathcal{U}$, and set
\begin{align}\label{eqn:control_prob}
    \Phi\big(t, \vec{z}(t)\big) = \inf_u J(t, \vec{z}, u), \quad \text{ s.t. } \eqref{eqn:dynamics}. 
\end{align}
Here, $\Phi$ is called the value function or the optimal cost-to-go. A solution $u^*$ incurring this minimal cost is called an optimal control, and the corresponding state variable $\vec{z}^*$ is called an optimal trajectory.

\begin{figure*}[htbp]
\floatconts
    {fig:local-sol}
    {\caption{\textbf{Local solution} approach drives initial \textbf{pathological} HH system, with $\Bar{g}_{\rm Na} = 380$, toward \textbf{normal} behavior by injecting control stimuli.}}
    {%
    \subfigure[Membrane Potential][c]{\label{fig:local_pathological_action_potential}%
        \includegraphics[width=0.3\linewidth, height=0.24\textwidth]{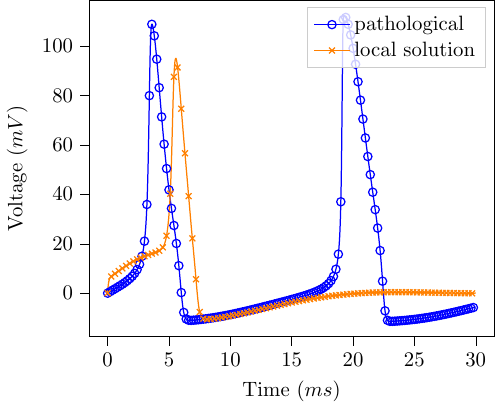}}
    \subfigure[Gating Variables][c]{\label{fig:local_pathological_gating_variables}%
        \includegraphics[width=0.3\linewidth, height=0.24\textwidth]{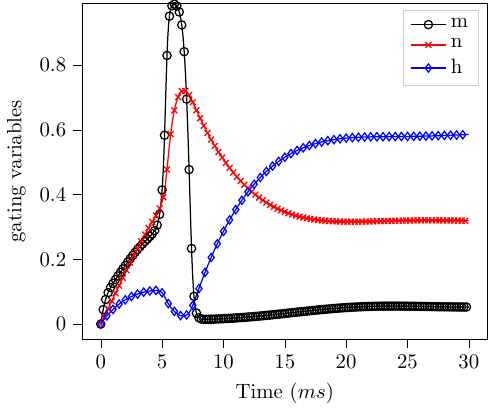}}
    \subfigure[Controls (stimuli)][c]{\label{fig:local_pathological_controls}%
        \includegraphics[width=0.3\linewidth, height=0.24\textwidth]{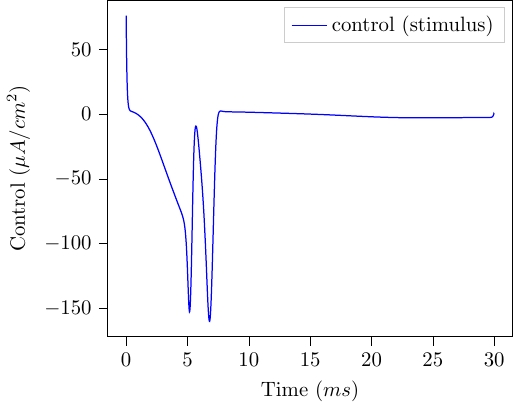}}
    }
\end{figure*}

\subsection{Optimal Control Solution Approaches}
\label{subsec:control_approaches}

Controls problems like \eqref{eqn:control_prob} can be solved using local or global optimization techniques. 

\paragraph{Local solution methods} aim to find an optimal control policy based on a fixed initial state. These methods are typically used for \textit{open-loop }control, which involves determining the control inputs (i.e., the optimal control trajectory) in advance, often through optimization techniques, and then applying these predetermined control inputs to the system without real-time feedback. However, re-computation or adaptation of the control strategy may be necessary when unexpected changes or perturbations occur in the system. This approach mirrors the conventional DBS practice of clinicians manually adjusting stimulation parameters to provide effective therapy for a patient. As the patient's condition evolves, additional appointments become necessary to fine-tune the treatment for optimal outcomes \citep{parastarfeizabadi2017advances, yu2020review}. 

A representative example of open-loop control is the \textit{all-at-once} Interior Point Method (IPM) \citep{nocedal2006interior}, a numerical optimization technique that determines the control inputs for the entire trajectory in advance, without incorporating real-time feedback from the system. We employ this approach as our baseline controller and iteratively solve the optimization problem \eqref{eqn:control_prob} by considering the entire control trajectory as a single optimization variable, simultaneously optimizing both the control input and the system's state trajectory.

We also consider Pontryagin's Maximum Principle (PMP) \citep{pontryagin2018mathematical}, a \textit{local solution} method that provides first-order necessary conditions for optimality and yields an open-loop optimal control law \citep{onken2021neural,li2022neural,flemingsoner06} as follows,
\begin{align}\label{eqn:pmp_control}
    u^*(s) \in \argmax_u \mathcal{H}\big(s,\vec{z}^*(s),\vec{p}(s), u(s)\big),
\end{align}
for every $s \in [0,T]$.
Here, for an adjoint state $\vec{p}$, the Hamiltonian $\mathcal{H}$ of the system is defined as
    \begin{align}\label{eqn:hamiltonian}
        \mathcal{H}(t,\vec{z},\vec{p}, u) &= - L(t, \vec{z}, u)- \vec{p}^\top \big[f(t,\vec{z})+\vec{e}_1 u\big].
    \end{align}
Solving this system yields a solution that is limited to a specific initial state, necessitating re-computation for different initial states or trajectory perturbations.

\paragraph{Global solution methods} solve the control problem for any given initial state, often yielding controls in \textit{closed-loop} or \textit{feedback} form. The control inputs are adjusted in real-time based on the system's current state and its response to previous control actions, enabling closed-loop control. This is relevant for critical applications such as DBS where real-time adaption to a patient's ongoing symptoms is crucial.

The Hamilton-Jacobi-Bellman (HJB) equation \citep{kirk2004optimal} (also known as the Dynamic Programming (DP) equation \citep{bellman1966dynamic}) is a global solution method used to derive the optimal control policy over a larger state space. Given any initial state, we can solve \eqref{eqn:control_prob} and produce an optimal control that leads to global convergence of the closed-loop system. DP states that the value function $\Phi$ satisfies the HJB equation \citep{flemingsoner06}, 
\begin{align}\label{eqn:hjb}
	\left\{ \begin{array}{cc}
		-\partial_t \Phi\big(t, \vec{z}(t)\big) + \mathcal{H}\big(t,\vec{z},\nabla_{\vec{z}} \Phi\big(t, \vec{z}(t)\big), u \big) = 0,\\
		\Phi\big(T, \vec{z}(T)\big) = G\big(\vec{z}(T)\big). 
	\end{array}\right.
\end{align}

In the context of DBS, the HJB equation guides the development of broad control policies that account for the entire state space. It provides a more comprehensive approach to optimizing stimulation across a wide range of potential scenarios.

\subsection{Semi-Global Approach}\label{subsec:semi_global}
While solving the HJB equation is generally fast and feasible when the state dimension $\leq 3$, it is prone to the CoD for higher dimensional problems. To mitigate this, we empoly a semi-global approach that uses neural networks to approximate the value function defined by the HJB equation, owing to their universal approximation properties \citep{han2018solving, kunisch2021semiglobal, onken2021neural, Exarchos-2018, pereira20a}.

This approach leverages the connection between the PMP and the HJB equation, which has long been established in control theory \citep{cernea2005connection,flemingsoner06} as 
\[ \vec{p}(t) = -\nabla_{\vec{z}} \Phi\big(t, \vec{z}(t)\big). \]
This relation, at optimality, helps obtain the optimal control $u^*(s)$ from the value function $\Phi$ at any given time via the feedback form, 
\begin{align}\label{eqn:feedback_form}
    u^*(s) \in \argmax_u \mathcal{H}\big(s,\vec{z}^*(s), \nabla_{\vec{z}} \Phi\big(s, \vec{z}^*(s)\big), u\big).
\end{align}
Assuming that we have a close-form solution for $u^*$ in \equationref{eqn:feedback_form} (which is true for our case), then the optimal controls can be recovered in real time using $\Phi$ and its gradient, $\nabla_{\vec{z}} \Phi$. This is ideal for clinical applications such as neurostimulation where swiftly computing controls for different times or states in real time is highly desirable. 

\paragraph{Ties to Reinforcement Learning} It is also important to note that the value function, $\Phi$, in the HJB equation is closely related to the value function estimated by methods such as actor-critic in Reinforcement Learning (RL) \citep{konda1999actor}. Ties between control theory and RL have been cemented by seminal works including \citep{bertsekas2019reinforcement, recht2018tour, sutton1992reinforcement}. In RL, the optimal policy maximizes $\Phi$, whereas in OC, we seek a control policy that minimizes $\Phi$. While scarce, existing applications of RL to neuroscience and deep brain stimulation include \citep{botvinick2020deep, gao2020model, Krylov_2020, Lu2020application}. These works demonstrate the potential of control-based techniques in developing innovative solutions for long-standing problems in healthcare and beyond.

\section{Learning Problem}
\label{sec:NN_approach}

Following \sectionref{subsec:semi_global}, we approximate the value function, $\Phi$, in \equationref{eqn:control_prob} using a neural network, ${\rm N}(\vec{y};\vec{\theta}_{\rm N})$, with parameters $\vec{\theta}$ as follows
\begin{align}\label{eqn:nn_approximation}
    \Phi_{\vec{\theta}}(\vec{y}) = \vec{w}^\top {\rm N}(\vec{y};\vec{\theta}_{\rm N}) + \frac{1}{2} \vec{y}^\top (\vec{A}^\top \vec{A})\vec{y}+\vec{b}^\top \vec{y}+c, \nonumber\\
    \text{where }\vec{\theta} = (\vec{w}, \vec{\theta}_{\rm N}, \vec{A}, \vec{b}, c),
\end{align}    
with space-time inputs $\vec{y} = \left(s, \vec{z}(s)\right) \in \mathbb{R}^{d+1}$. Here, $d$ is the state dimension ($=4$) and $\vec{\theta}$ consists of all trainable weights: $\vec{w} \in \mathbb{R}^m, \vec{\theta}_{\rm N} \in \mathbb{R}^{p_{\rm N}}, \vec{A} \in \mathbb{R}^{d \times (d+1)}, \vec{b} \in \mathbb{R}^{d+1}, \text{ and } c \in \mathbb{R}$, where $\vec{A}, \vec{b}$ and $c$ model linear dynamics and neural network ${\rm N}(\vec{y};\vec{\theta}_{\rm N}): \mathbb{R}^{d+1} \rightarrow \mathbb{R}^m$ models nonlinear dynamics. $p_{\rm N}$ denotes the number of parameters of the NN.

\begin{figure*}[htbp]
\floatconts
    {fig:semiglobal-sol}
    {\caption{\textbf{Semi-global solution} approach leveraging neural networks successfully drives initial \textbf{pathological} HH system, with $\Bar{g}_{\rm Na} = 380$, toward \textbf{normal} behavior by injecting control stimuli.}}
    {%
    \subfigure[Membrane Potential][c]{\label{fig:NN_pathological_action_potential}%
        \includegraphics[width=0.3\linewidth, height=0.24\textwidth]{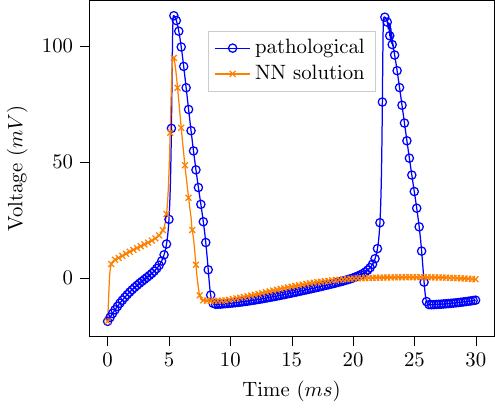}}
    \subfigure[Gating Variables][c]{\label{fig:NN_pathological_gating_variables}%
        \includegraphics[width=0.3\linewidth, height=0.24\textwidth]{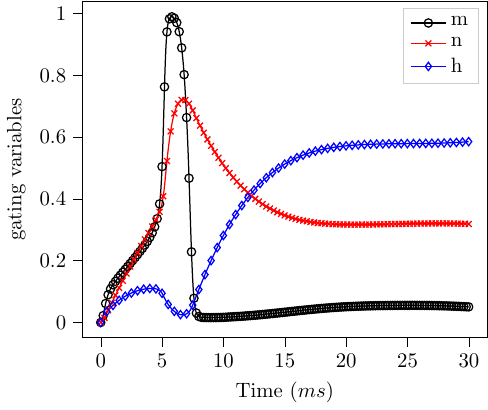}}
    \subfigure[Controls (stimuli)][c]{\label{fig:NN_pathological_controls}%
        \includegraphics[width=0.3\linewidth, height=0.24\textwidth]{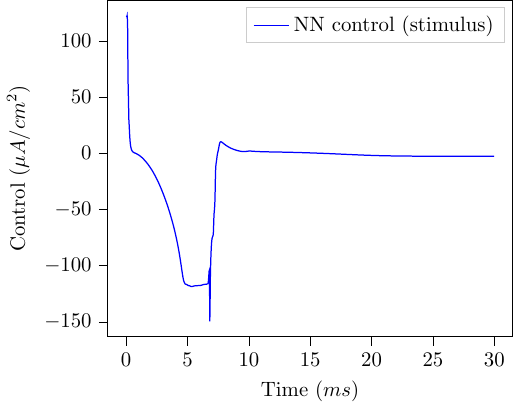}}
    }
\end{figure*}

To learn the parameters $\vec{\theta}$ of the NN, we first sample some initial states $\vec{x}$ from the Gaussian distribution, $\rho$ with mean $0$ and a variance of $10$, and approximately solve the following minimization problem
\begin{align}\label{eqn:learning_problem} 
    \min_{\vec{\theta}}  \mathbb{E}_{z_0 \sim \rho}\left\{\begin{multlined} \ell(T) + g\big(\vec{z}(T)\big) + \gamma_1 c_{\rm HJB}(T) \\+ \gamma_2 \big|\Phi_{\vec{\theta}} \big(T, \vec{z}(T)\big) - g\big(\vec{z}(T)\big)\big| \end{multlined} \right\},
\end{align}   
subject to 
\begin{align}
    \begin{split}
        \partial_s &\begin{pmatrix}
            \vec{z}(s) \\ 
            \ell(s) \\ 
            c_{\rm HJB}(s)
        \end{pmatrix} 
        =\\
        &\begin{pmatrix}
            -\nabla_{\vec{p}}\mathcal{H}\big(s, \vec{z}(s), \nabla_{\vec{z}} \Phi_{\vec{\theta}}\big(s, \vec{z}(s)\big), u^* \big) \\
            L_{\vec{z}}(s)\\
            R_{\vec{z}}(s)
        \end{pmatrix},
    \end{split}
\end{align}
for $s\in (0,T)$, initialized with $\ell(0) = c_{\rm HJB}(0) = 0$ and $\vec{z}(0) = \vec{x}$. Here,
\begin{align*}
 L_{\vec{z}}(s)&=\\ &\nabla \Phi_{\vec{\theta}} \big(s, \vec{z}(s) \big) \cdot \nabla_{\vec{p}} \mathcal{H} \big( s, \vec{z}(s), \nabla \Phi_{\vec{\theta}}\big(s, \vec{z}(s)\big), u^* \big),
 \end{align*}
 and
 \begin{align*}
R_{\vec{z}}(s) = \big| -\partial_t \Phi_{\vec{\theta}}& \big(s, \vec{z}(s)\big) + \\  &\mathcal{H}\big(s,\vec{z}(s),\nabla_{\vec{z}} \Phi_{\vec{\theta}}\big(s, \vec{z}(s)\big), u^*\big) \big|,
\end{align*}
with $\ell$ denoting the accumulated running cost.

The hyperparameters $\gamma_1, \gamma_2 \geq 0$ balance minimization of the control objective and HJB penalization.

\section{Numerical Experiments}
\label{sec:exp}

In this section, we examine two stimulation strategies: a baseline controller that applies the all-at-once IPM to solve the control problem in an open-loop manner, i.e., with no system feedback, and a controller that leverages neural networks for optimized controls over a broad state space (semi-global method). Our NN-based controller architecture has depth $2$, width $64$, a learning rate of $0.005$, and uses the ADAM optimizer. We evaluate robustness by introducing minor perturbations or shocks to the HH system, \equationref{eqn:dynamics}.

\subsection{Neuromodulatory Effects}
We simulate the HH model for both normal and pathological conditions, based on parameters in \tableref{table:HH_params}. The former is achieved by modeling with the original parameters defined in \citep{hodgkin1952quantitative}, and the latter is obtained by varying the parameters, specifically targeting $\Bar{g}_{\rm Na}$. This distorts the ion flow across the neural membrane, which in turn disrupts the action potential generation, see \figureref{fig:normal_vs_pathological_action_potential}. In this scenario, the objective is to develop control strategies that restore the disrupted system to its normal condition while minimizing the energy required for input current injection. This is valuable in countering pathological neural activity in conditions like Parkinson's disease or chronic pain \citep{little2013adaptive, shirvalkar2018closed}.

For both approaches, we initialize the simulation by perturbing the neuron's resting state (adding Gaussian noise scaled by a factor of $10$ to $V_{\rm m} = 0$ while setting the other state variables to $0$). This captures the effect of the aforementioned pathological activity. \figureref{fig:local-sol} shows the state trajectories and controls learned by the baseline while driving the pathological states toward the normal states. Likewise, \figureref{fig:semiglobal-sol} shows the state trajectories and controls learned by the neural network approach on the same problem. Both approaches restore the pathological HH system to normal conditions, albeit with varying stimuli.

\begin{table}[htbp]
    \floatconts
    {table:both_objF}
    {\caption{Running and terminal costs for single instance shown in \figureref{fig:local-sol,fig:semiglobal-sol}.}}
    {
        \begin{small}
        \begin{sc}
        \begin{tabular}{lccccr}
        \toprule
         & running ($\ell$) & terminal ($G$) & total \\
        \midrule
        $\text{IPM}$ &  45725.6861 & 0.0279 & 45725.7140 \\
         $\text{NN}$ & 46890.1119  & 0.0088 & 46890.1207 \\
        \bottomrule
        \end{tabular}
        \end{sc}
        \end{small}
    }
\end{table}

\subsection{Suboptimality} \label{sec:baseline}
The baseline IPM controller modulates neuronal activity for a fixed initial state $x$ at a time. We consider the solution obtained using this approach as the ground-truth optimal solution. To compare the performance of the semi-global NN approach, we evaluate its solution for the initial state and compute its suboptimality relative to the baseline solution. \tableref{table:both_objF} compares the objective values of both methods on $\vec{x} = [0, 0, 0, 0]$. The NN-based controller achieves near optimal performance, considering it solves the control problem for a larger state space compared to the baseline.

\begin{figure}[htbp]
\floatconts
{fig:sub_opt}
{\caption{Suboptimality comparison between NN approach and local solution method for $\xi$ in $[-40, 40]$. \textbf{Left}: Overall trend. \textbf{Right}: Zoomed-in view on trained sample points.}}
{
    \subfigure[Overall trend][c]
    {
        \label{fig:overall_subopt}%
        \includegraphics[width=0.425\linewidth]{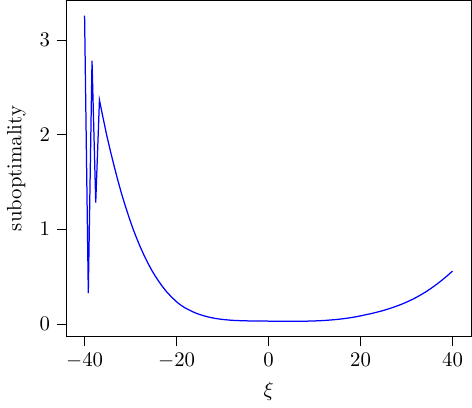}%
    }
    \qquad
    \subfigure[Trained Sample Points][c]
    {
        \label{fig:zoomed_subopt}%
        \includegraphics[width=0.425\linewidth]{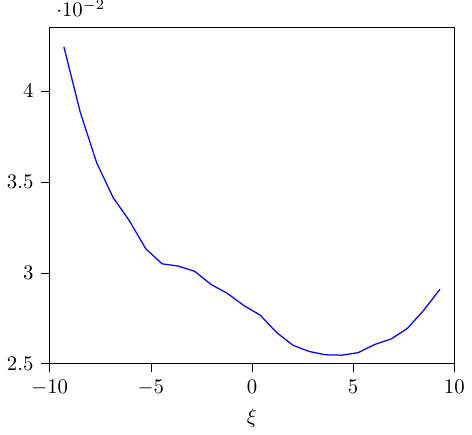}
    }
}
\end{figure}

We further evaluate the performance of the NN solution by measuring its suboptimality across various initial conditions $(\xi,0,0,0)$. To achieve this, we uniformly sample $100$ values of $\xi$ between $-40.0$ and $40.0$, representing the membrane potential $V_{\rm m}$. For each $(\xi,0,0,0)$ combination, where the local solution is considered optimal for every initial state, we calculate the suboptimality of the NN solution relative to the local solution. It is worth noting that our network is specifically trained for $\xi$ values within the range of $[-10, 10]$. Therefore, evaluating the performance across a broader range of $\xi$ values also provides insights into how our model performs on out-of-distribution samples. The results of this evaluation are depicted in \figureref{fig:sub_opt}.
\begin{figure}[htbp]
    \floatconts
    {fig:shock}
    {\caption{The NN handles a shock to the system dynamics and recovers the optimal trajectory.}}
    {
    \subfigure[$V_{\rm m}$ with shock][c]
    {
        \label{fig:Vm_shock}%
        \includegraphics[width=0.6\linewidth]{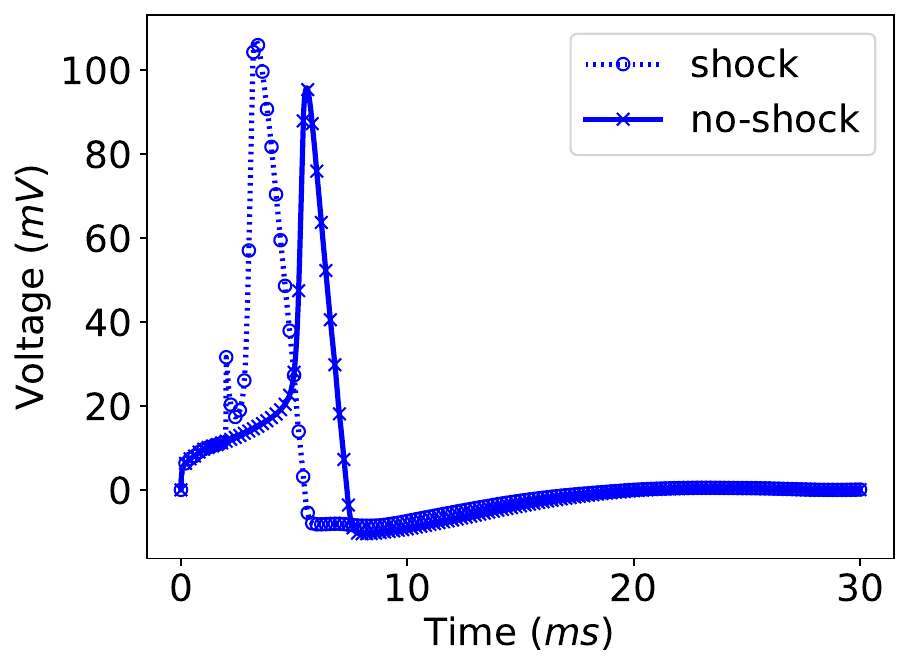}%
    }  
    \\
    \subfigure[$m, n, h$ with shock][c]
    {
        \label{fig:gating_vars_shock}%
        \includegraphics[width=0.6\linewidth]{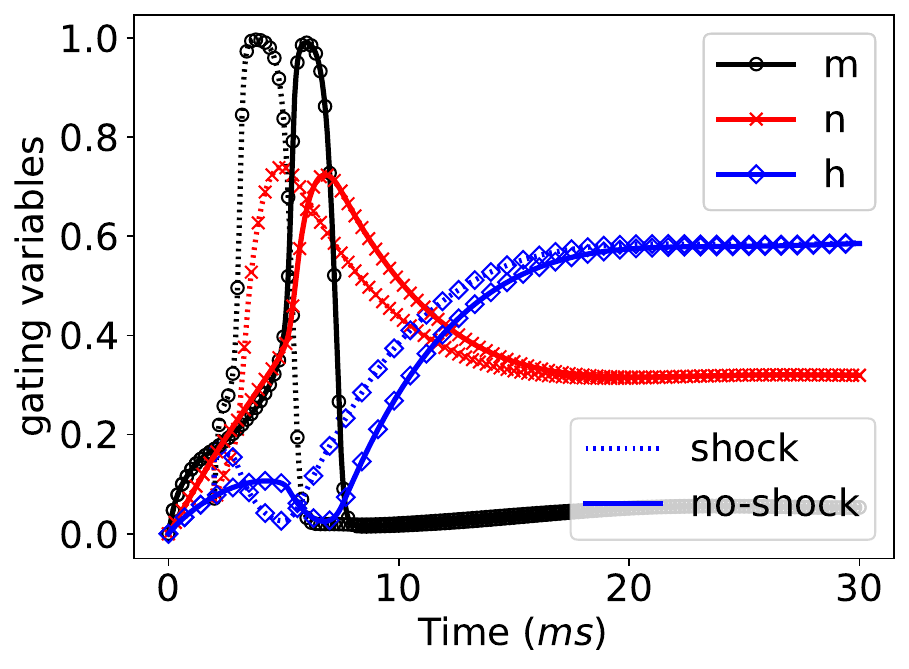}
    } 
    }
\end{figure}
\subsection{Robustness to shocks}
We investigate the robustness of the NN control policy to shocks or disturbances in the system dynamics. In \figureref{fig:shock}, we illustrate two trajectories: a normal trajectory and one with an added shock, both starting from the initial state $\vec{x} = [0, 0, 0, 0]$. Remarkably, the NN is capable of recovering the optimal trajectory even after the shock. This can be attributed to the semi-global approximation of the value function employed by the NN. Additionally, since the NN is trained offline, it can handle real-time disturbances without the need for re-computation, unlike the baseline approach.

\section{Conclusion}
This work explores the intersection of neuronal dynamics, optimal control theory, and machine learning. We focus on controlling neuronal dynamics characterized by the Hodgkin-Huxley model, a nonlinear system of differential equations that provides a detailed and complex representation of the spiking behavior of neurons. This has direct implications for biomedical applications like deep brain stimulation, where modulating neuronal dynamics to achieve therapeutic benefit with minimal energy expenditure is desired. Following existing approaches \citep{han2018solving, kunisch2021semiglobal, onken2021neural}, we approximate the value function with a neural network and combine the PMP and HJB approaches to obtain feedback control policies. This allows for real-time computations of controls based on the current state and leads to overall improved robustness to disturbances, as demonstrated in our numerical experiments. This signifies a promising research direction in closed-loop DBS, for enabling a dynamic and personalized approach to treatment. By combining optimal control theory and machine learning, we highlighted new possibilities for interdisciplinary research, particularly toward the design of adaptive neurotechnology and closed-loop systems. In future work, we expect to extend it to multi-compartment models featuring populations of neurons. These models comprise hundreds of neurons interacting across various neural structures, which pose exciting computational challenges that remain to be fully explored in these applications.

\acks{This material is based on work supported by several grants including the 2021 Google PhD Fellowship in Computational Neural and Cognitive Sciences, AFOSR grant FA9550-20-1-0372, US DOE Office of Advanced Scientific Computing Research Field Work Proposal 20-023231, and NSF awards DMS 1751636 and DMS 2038118.}

\bibliography{ml4h_references}

\appendix

\section{Background: Optimal Control theory}\label{apd:oc_background}
We present a brief background to Optimal Control (OC) theory relevant to this work. For a complete introduction to the topic, please refer to canonical texts such as \citet{evans1983introduction, kirk2004optimal,yong1999stochastic}.
 
 OC is branch of applied mathematics that aims to find a control strategy or policy able to manipulate a dynamical system optimally over time in order to achieve a desired target state(s) of the system. OC theory has applications in various fields, including engineering, economics, and healthcare. Solving a control problem generally requires \textbf{minimizing a cost function} (performance criterion) associated with the dynamical system. Mathematically, given a dynamical system, defined over a fixed finite time-horizon $T$,
\begin{align}\label{eqn:apd_dynamics}
\begin{split}
\frac{d \vec{z}}{dt}(t) &= f\big(t, \vec{z}(t), \vec{u}(t)\big),\; \quad 0 \leq t \leq T\\
\vec{z}(0) &= \vec{x},
\end{split}
\end{align}
where $\vec{z}(t) \in \mathbb{R}^d$ represents the \textbf{state} of the system with dimension $d$ at time $t$, with $\vec{x}$ being the initial state of the system. $\vec{u}(t): [0, T] \to \mathbb{R}$ represents the external \textbf{control} provided as input by a controller (agent) at time $t$. The function $f: [0,T] \times \mathbb{R}^d \to \mathbb{R}^d$ governs the evolution of the state dynamics and is assumed to be \textit{known}. 

The primary objective of the OC problem is to find the control that drives the dynamical system toward a reference/target state, $\vec{z}^*$, while incurring the minimum possible cost. This involves specifying an \textbf{objective functional} (cost/utility function), $J$, as follows
\begin{align}\label{eqn:apd_cost_function}
    J(t, \vec{z}, \vec{u}) &= G\big(\vec{z}(T)\big) + \int_t^T L\big(s, \vec{z}(s), \vec{u}(s)\big) ds,
\end{align}
where $L: [0, T] \times \mathbb{R}^d \times \mathbb{R} \to \mathbb{R}$ is the \textit{Lagrangian} or \textit{running} cost and $G: \mathbb{R}^d \to \mathbb{R}$ is the \textit{terminal} cost. $L$ accumulates the cost of controlling the system and expending energy at each time step. $G$ penalizes the distance between the final state $\vec{z}(T)$ and the given target state of the system, $\vec{z}^*(T)$, at final time. 

We are interested in finding a control(s) that incurs the minimal cost over all admissible controls $\vec{u} \in \mathcal{U}$,
\begin{align}\label{eqn:apd_control_prob}
    \Phi\big(t, \vec{z}(t)\big) = \inf_{\vec{u}} J(t, \vec{z}, \vec{u}), \quad \text{ s.t. } \eqref{eqn:apd_dynamics}.
\end{align}
 A solution $\vec{u}^*$ incurring this minimal cost is called an optimal control, and the corresponding state variable $\vec{z}^*$ is called an optimal trajectory.

\subsection{All-at-once Interior Point Method}\label{apd:baseline}
The All-at-Once Interior Point Method (IPM) \citep{nocedal2006interior} aims to find the optimal control $u^*(t)$ by solving the following optimization problem:

\begin{align}
\begin{split}
& \min_{\vec{u}(t), \vec{z}(t)} J\left(t, \vec{z}, \vec{u}\right) \\
& \text{s.t. } \\
& \quad \frac{d \vec{z}}{dt}(t) = f\big(t, \vec{z}(t), \vec{u}(t)\big)\\
& \quad g(\vec{z}(t), \vec{u}(t)) + s = 0 \\
& \quad h(\vec{z}(T)) = 0 \\
& \quad s \geq 0
\end{split}
\end{align}
This formulation introduces slack variables $s$ to handle inequality constraints, and the barrier function is often used to incorporate them into the objective function. The IPM is solved iteratively, where at each iteration, the first-order necessary conditions are considered and barrier parameters are updated \citep{nocedal2006interior}.

\end{document}